\documentclass[11pt,twocolumn]{IEEEtran}

\pdfoutput=1  
\usepackage{cite,url}

\usepackage[usenames,dvipsnames]{color}
\usepackage{graphicx,wrapfig,subfigure}

\usepackage[cmex10]{amsmath}
\usepackage{amsfonts,amssymb,mathrsfs}   
\usepackage[colorlinks,bookmarksopen,bookmarksnumbered,citecolor=red,urlcolor=red]{hyperref}
\usepackage{soul}

\newlength{\noteWidth}
\long\def\notes#1{\ifinner
           {\footnotesize #1}
           \else
           \marginpar{\parbox[t]{\noteWidth}{\raggedright\footnotesize #1}}
       \fi\typeout{#1}}

\def\notes#1{\typeout{read notes: #1}}  

\def\spm#1{\notes{SPM:  #1}}

\newcommand*{\qed}{\nobreak\hfill\ensuremath{\square}}

\newcounter{rmnum}

\newenvironment{romannum}{\begin{list}{{\upshape (\roman{rmnum})}}{\usecounter{rmnum}
\setlength{\leftmargin}{8pt}
\setlength{\rightmargin}{8pt}
\setlength{\itemsep}{2pt}
\setlength{\itemindent}{-1pt}
}}{\end{list}}
 
\newcounter{anum}

\def\poleshift{\varrho}

\def\tily{\tilde{y}}

\def\util{\mathchoice{\mbox{\small$\cal U$}}%
{\mbox{\small$\cal U$}}%
{\mbox{$\scriptstyle\cal U$}}%
{\mbox{$\scriptscriptstyle\cal U$}}}

\def\bfgamma{\bfmath{\gamma}}

\def\Cov{\text{Cov}\,}   
\def\diag{\text{diag}\,}   

\def\Spx{\textsf{S}}
\def\Ebox#1#2{%
\begin{center}
\includegraphics[width= #1\hsize]{figures/#2} 
\end{center}}

\def\bfmtily{\widetilde{\bfmath{y}}}

\def\Fig#1{Fig.~\ref{#1}}

\def\ind{\field{I}}

\def\Re{\field{R}}

\def\piload{\Gamma}
\def\Health{{\cal L}}
\def\health{\ell}
\def\psd{\text{S}}

\def\bfmath#1{{\mathchoice{\mbox{\boldmath$#1$}}%
{\mbox{\boldmath$#1$}}%
{\mbox{\boldmath$\scriptstyle#1$}}%
{\mbox{\boldmath$\scriptscriptstyle#1$}}}}

\def\bfmr{\bfmath{r}}

\def\bfmy{\bfmath{y}}

\def\bfmw{\bfmath{w}}

\def\bfmB{\bfmath{B}}

\def\bfmD{\bfmath{D}}

\def\bfmX{\bfmath{X}}

\def\bfmY{\bfmath{Y}}

\def\bfmhhaY{\bfmath{\hhaY}} 
\def\bfmhhaY{\hbox to 0pt{$\widehat{\bfmY}$\hss}\widehat{\phantom{\raise 1.25pt\hbox{$\bfmY$}}}}

\def\bfPhi{\bfmath{\Phi}}

\def\bfzeta{\bfmath{\zeta}}

\newtheorem{theorem}{Theorem}[section]

\newtheorem{proposition}[theorem]{Proposition}
\newtheorem{lemma}[theorem]{Lemma}

\def\Proposition#1{Proposition~\ref{#1}}

\def\Section#1{Section~\ref{#1}}

\def\state{{\sf X}}

\newcommand{\field}[1]{\mathbb{#1}}

\def\Co{\field{C}}

\def\Re{\field{R}}

\def\nat{\field{Z}_+}

\def\bary{{\overline {y}}}

\def\Prob{{\sf P}}

\def\Expect{{\sf E}}

\def\eqdef{\mathbin{:=}}

\def\transpose{{\hbox{\it\tiny T}}}

\def\clE{{\cal E}}

\title{Individual Risk in Mean Field Control 
\\
with
Application to Automated Demand Response}

\author{Yue Chen, 
Ana Bu\v{s}i\'c, and Sean Meyn
\thanks{This research is supported by the NSF grants CPS-0931416 and CPS-1259040, the French National Research Agency grant ANR-12-MONU-0019, and US-Israel BSF Grant 2011506.}
\thanks{
Y.C. and S.M. are with the Department of Electrical and Computer
Engg.\ at the University of Florida, Gainesville. A.B.\ is with Inria and the Computer Science Dept. of \'Ecole Normale Sup\'erieure, Paris, France.}%
}

\begin{document}

\maketitle

\spm{More refs?  Roland sent this:
kizmal13
Mean Field Based Control of Power System Dispersed Energy Storage
Devices for Peak Load Relief
Arman C. Kizilkale Roland P. Malham´e}

\setcounter{page}{0}
\thispagestyle{empty}

\begin{abstract} 

Flexibility of energy consumption can be harnessed for the purposes of ancillary services in a large power grid.

In prior work by the authors a randomized control architecture is introduced for individual loads for this purpose. 
In examples it is shown that the control architecture can be designed so that control of the loads is easy at the grid level: Tracking of a balancing authority reference signal is possible, while ensuring that the quality of service (QoS) for each load is acceptable \textit{on average}. The analysis was based on a mean field limit (as the number of loads approaches infinity), combined with an LTI-system approximation of the aggregate nonlinear model. 


This paper examines in depth the issue of individual risk in these systems. 
The main contributions of the paper are of two kinds:

\textit{Risk is modeled and quantified}
\begin{romannum}
\item[(i)]
The average performance is not an adequate measure of success. It is found empirically that a histogram of QoS is approximately Gaussian, and consequently each load will eventually receive poor service.

\item[(ii)]
The variance can be estimated from a refinement of the LTI model that includes a white-noise disturbance; 
variance is a function of the randomized policy, as well as the power spectral density of the reference signal.
\end{romannum}

\textit{Additional local control can eliminate risk}
\begin{romannum}

\item[(iii)]
The histogram of QoS is \textit{truncated} through this local control, so that strict bounds on service quality are guaranteed. 

\item[(iv)]
This has insignificant impact on the grid-level performance, beyond a modest reduction in capacity of ancillary service. 
\end{romannum}

\end{abstract}


\section{Introduction} 
\label{s:intro}

The power grid requires regulation to ensure that supply matches demand.  Regulation is required by each
 \textit{balancing authority} (BA) on multiple time-scales, corresponding to the time-scales of volatility of both supply and demand for power.  Resources that supply these regulation services are collectively known as  \textit{ancillary services}.  FERC orders 755 and 745 are intended to provide incentives for the provision of these services.  
The need for regulation resources increases significantly with greater penetration of renewable generation from wind or sun; 
this is observed today in regions of the world with significant generation from these sources.

A number of papers have explored the potential for extracting ancillary service through the inherent flexibility of loads.
Examples of loads with sufficient flexibility to provide service to the grid are aluminum manufacturing,
data centers,  plug-in electric vehicles,  heating and ventilation (HVAC), and water pumping for irrigation 
\cite{malcho85,matkoccal13,kizmal13,meybarbusyueehr14,haolinkowbarmey14}. 
\spm{Reference IBM in journal version}

An aggregate of loads can be modeled as a dynamical system.  The dynamics of the aggregate depends upon the physical properties of the loads, as well as the mechanism that is used to engage them.     For example,  with the use of price signals, load response may exhibit  threshold behavior in which most loads shut down when the price exceeds some value.   Even with direct load control, there may be delay and dynamics, so harnessing ancillary services from flexible loads amounts to a control problem: The BA wishes to design some signal to be broadcast to loads, based on measurements of aggregate power output, so that deviation in power consumption tracks a reference signal -- See \cite{meybarbusyueehr14} for detailed discussion.

For many types of loads it has been argued that a randomized control architecture at each load
can simplify this control problem    \cite{matkoccal13,meybarbusyueehr14}.  
\begin{figure}[h]
\Ebox{.95}{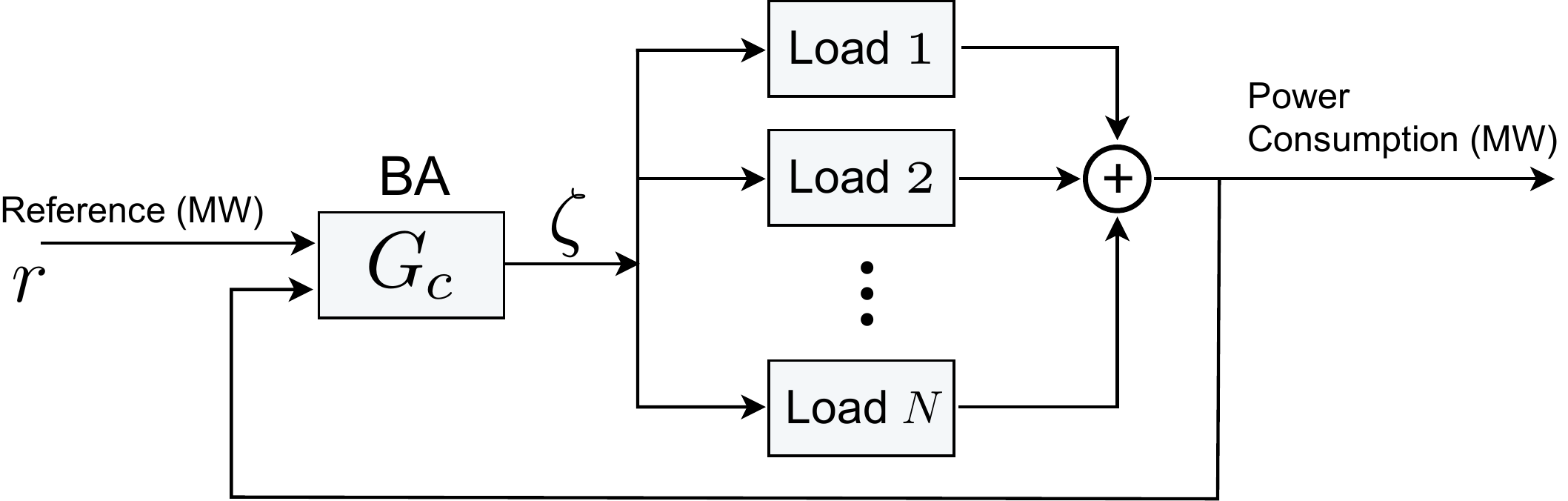} 
\vspace{-.15cm}
\caption{The control architecture: command $\bfzeta$ is computed at a BA,
 and transmitted to each load. The control decision at a load is based only on its own state and the signal $\bfzeta$.}
\vspace{-.25cm}
\label{fig:arch} 
\end{figure} 

\Fig{fig:arch} shows a schematic of the control architecture adopted in  \cite{meybarbusyueehr14}, 
in which each load operates according to a randomized policy based on its internal state, and a common control signal $\bfzeta$. 
Theoretical results and examples in this prior work demonstrate that local randomized policies can be designed so that control of the loads is easy at the grid level.     The analysis was based on a mean field limit (as the number of loads approaches infinity), combined with a linear time-invariant (LTI) system approximation of the aggregate nonlinear model.   In the examples considered in this prior work, the linear approximation was found to be minimum phase, which is why simple linear error feedback could be applied to achieve nearly perfect tracking.

Absent in prior work is any detailed analysis of risk for an individual load.  In the setting of  \cite{meybarbusyueehr14} it can be argued that  the quality of service (QoS) for each load is acceptable \textit{on average}, but the volatility of QoS has not been addressed to-date.  Strict bounds on QoS are addressed in \cite{leboudec11} for a deterministic model.  The \textit{service curves} considered there are of similar flavor to the QoS metrics used in the present work.  
\spm{NEW}


\begin{figure*} 
\vspace{.2cm}
\Ebox{.9}{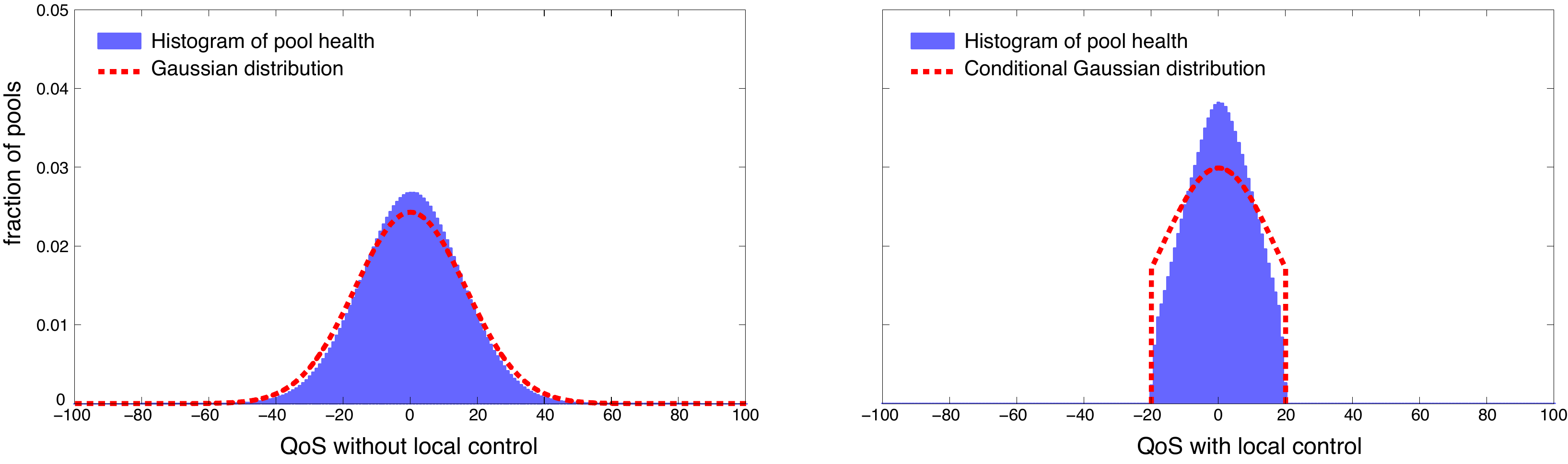}
\vspace{-.2cm}
\caption{Histogram of the discounted QoS \eqref{e:Health}: with and without local opt-out control.}
\label{fig:BothHist} 
\end{figure*} 

This paper examines in depth the issue of individual risk. 
The main contributions of the paper are of two kinds.

In the first part of the paper we propose approaches to quantify QoS,  and methods to approximate its mean and variance so that we can quantify risk.   The variance is estimated based on a stochastic LTI model that approximates the dynamics of the time-dependent QoS function.  The additive disturbance in the LTI model is a function of the randomized policy used at each load and also the exogenous reference signal.  Estimation of the variance of QoS is thus reduced to an estimation of the  power spectral density of these disturbances.  These theoretical results are developed in \Section{s:MFMload}.

A simple approach is proposed to restrict QoS to pre-specified bounds:  A load will opt-out of service to the grid temporarily, whenever its QoS is about to exit pre-specified bounds.   This essentially eliminates risk since the histogram of QoS is restricted to these bounds.

Numerical results surveyed in \Section{s:num} confirm that the histogram of QoS is truncated through this local control, so that strict bounds on service quality are guaranteed.  This has \textit{insignificant} impact on the grid-level performance, beyond a modest reduction in capacity of ancillary service.

\Fig{fig:BothHist} shows a histogram of QoS based on simulation experiments described in \Section{s:num}.  The plot on the left hand side shows that the Gaussian approximation is a good fit with empirical results when there is no local opt-out control.  The figure on the right shows how the histogram is truncated when opt-out control is in place.

We begin with a brief survey of a portion of results from \cite{meybarbusyueehr14},  and a precise definition of QoS for a load. 

\section{Randomized control and mean-field models} 
\label{s:mfm}

\subsection{Randomized control}

The system architecture considered in this paper is
illustrated in \Fig{fig:arch}, based on the following components: 
\begin{romannum}
\item 
There are $N$ homogeneous loads that receive a common scalar command signal from the balancing authority, or BA, denoted $\bfzeta=\{\zeta_t\}$ in the figure.

\item
Each load evolves as a controlled Markov chain on the finite state space $\state=\{x^1,\dots,x^d\}$.
Its transition probability is determined by its own state, and the BA signal $\bfzeta$. The common dynamics are defined by a controlled transition matrix $\{P_\zeta : \zeta\in\Re\}$. For the $i$th load, there is a state process $\bfmX^i$ whose transition probability is given by,
\begin{equation}
\Prob\{X^i_{\tau+1} = x' \mid X^i_\tau = x ,\, \zeta_\tau=\zeta\} = P_{\zeta}(x,x') 
\label{e:Pzeta}
\end{equation}
for each $x, x'\in\state$.

\item
The BA has measurements of the other two scalar signals shown in the figure: The normalized aggregate power consumption $\bfmy$ and desired deviation in power consumption $\bfmr$, which is also normalized.

\end{romannum}
An approach to construction of $\{P_\zeta:\zeta\in\Re\}$ was proposed in  \cite{meybarbusyueehr14} based on information-theoretic arguments.  In the present paper we do not require a specific construction, but assume that $P_\zeta$ is continuously differentiable in the parameter $\zeta$.

The introduction of the time index $t$ in (i) and $\tau$ in (ii) is motivated by the possibility that the sampling time at the grid level may be faster than the sampling time at each load.   This is explained in \Section{s:ss}.

It is assumed that the power consumption at time $t$ from load $i$ is equal to some function of the state, denoted $\util(X^i_t)$.   The normalized power consumption is denoted,
\[
y_t^N = \frac{1}{N}\sum_{i=1}^N \util(X^i_t).
\]
The superscript is dropped unless dependency on $N$ must be emphasized.

The nominal behavior of each load is defined as the dynamics with $\bfzeta\equiv 0$.  In this case the loads are independent Markov chains.  It is assumed that the Markov chain is ergodic: It is uni-chain and aperiodic, so that $P_0$ has unique invariant distribution $\pi_0$.
The value $\bary^0\eqdef \sum_x \pi_0(x)\util(x)$ is the average nominal power usage.  On combining the ergodic theorem for Markov chains with the Law of Large Numbers for i.i.d.\ sequences we can conclude that $y^N_t\approx \bary^0$ when both $N$ and $t$ are large,  provided $\bfzeta\equiv 0$.  
 
It is assumed that the regulation signal is also normalized so that tracking amounts to choosing the signal $\bfzeta$ so that $r_t\approx \tily_t$ for all $t$, where $\tily_t = y_t - \bary^0$ is the deviation from nominal behavior. For example, we might use error feedback of the form, 
\begin{equation}
\zeta_t = G_c e_t,\qquad e_t = r_t-\tily_t   \, ,  
\label{e:zetaGc}
\end{equation} 
where $G_c$ is a transfer function.   The design of $G_c$ depends on a linear systems approximation for 
dependency of $\bfmy$ on $\bfzeta$.  This can be obtained by the construction of a mean-field model, as in the prior work 
\cite{coupertemdeb12,macalhis10,meybarbusyueehr14}.

\subsection{Mean-field model}

The mean-field model is based on consideration of the empirical distributions,
\begin{equation}
\mu^N_t(x)\eqdef \frac{1}{N}\sum_{i=1}^N  \ind\{ X^i_t =x \} ,\quad x\in\state.
\label{e:empDist}
\end{equation}
Under very general conditions on the input sequence $\bfzeta$,  it can be shown that the empirical distributions converge to a solution to the nonlinear state space model equations, 
\begin{equation}
\mu_{t+1} = \mu_t P_{\zeta_t}.
\label{e:MFM}
\end{equation} 
The interpretation remans the same:   $\mu_t$   represents the fraction of loads in each state, and sums to $1$.  
The output is denoted $y_t=\sum_x \mu_t(x)\util(x)$, which is the limit of the $\{y_t^N \}$ as $N\to\infty$.  See
 \cite{meybarbusyueehr14} for details.

The unique equilibrium with $\bfzeta\equiv 0$ is $\mu_t\equiv \pi_0$ and $y_t\equiv \bary^0 $.   Its linearization around this equilibrium is given by the linear state space model,
\begin{equation}
\begin{aligned}
 \Phi_{t+1} &= A \Phi_t + B \zeta_t
 \\
\gamma_t &= C \Phi_t
\end{aligned}
\label{e:LSSmfg}
\end{equation}
where
$A=P^\transpose_0$, $C$ is a row vector of dimension $d=|\state|$ with $C_i= \util(x^i) $ for each $i$, and $B$ is a $d$-dimensional column vector with entries $B_j = \sum_x\pi_0(x) \clE(x,x^j) $, where
the matrix $\clE$ is equal to the derivative,
\begin{equation}
\clE=\frac{d}{d\zeta} P_\zeta \Big|_{\zeta=0}.
\label{e:Pder}
\end{equation}

In the state equations \eqref{e:LSSmfg},  the state  $\Phi_t$ is $d$-dimensional,  and $\Phi_t(i)$ is intended to approximate $\mu_t(x^i) -\pi_0(x^i)$ for $1\le i\le d$.   The output $\gamma_t$ is an approximation of $\tily_t$.

\subsection{Super-sampling}
\label{s:ss}

An extension of this model is required in some applications to reduce delay at the grid level.  An approach called \textit{super sampling} is introduced without discussion in \cite{meybarbusyueehr14}.  Here we give a full description of this approach since it is required in the results that follow.   

We let $T$ denote the sampling time for each load, in units of minutes.  
We fix an integer $m>1$, and assume that $T_g\eqdef T/m$ is the sampling interval at the grid level.  The signal $\bfzeta$ is broadcast from the BA every $T_g$ minutes.   In applications,  $T$ might be equal to 30 minutes, while $T_g$ one minute.  

The $N$ loads are assumed to be equally divided into $m$  classes.
A load is said to be of class $k$,  where $0\le k\le m-1$,  
if the reference signal is sampled at times $\{iT +  kT_g : i\ge 0\}$.  
At grid-level time, the state of a load in class $k$ is constant on intervals of length $m$, 
with potential jumps only at the sampling times $\{iT +  kT_g : i\ge 1\}$.

\begin{figure}[h]
\Ebox{.95}{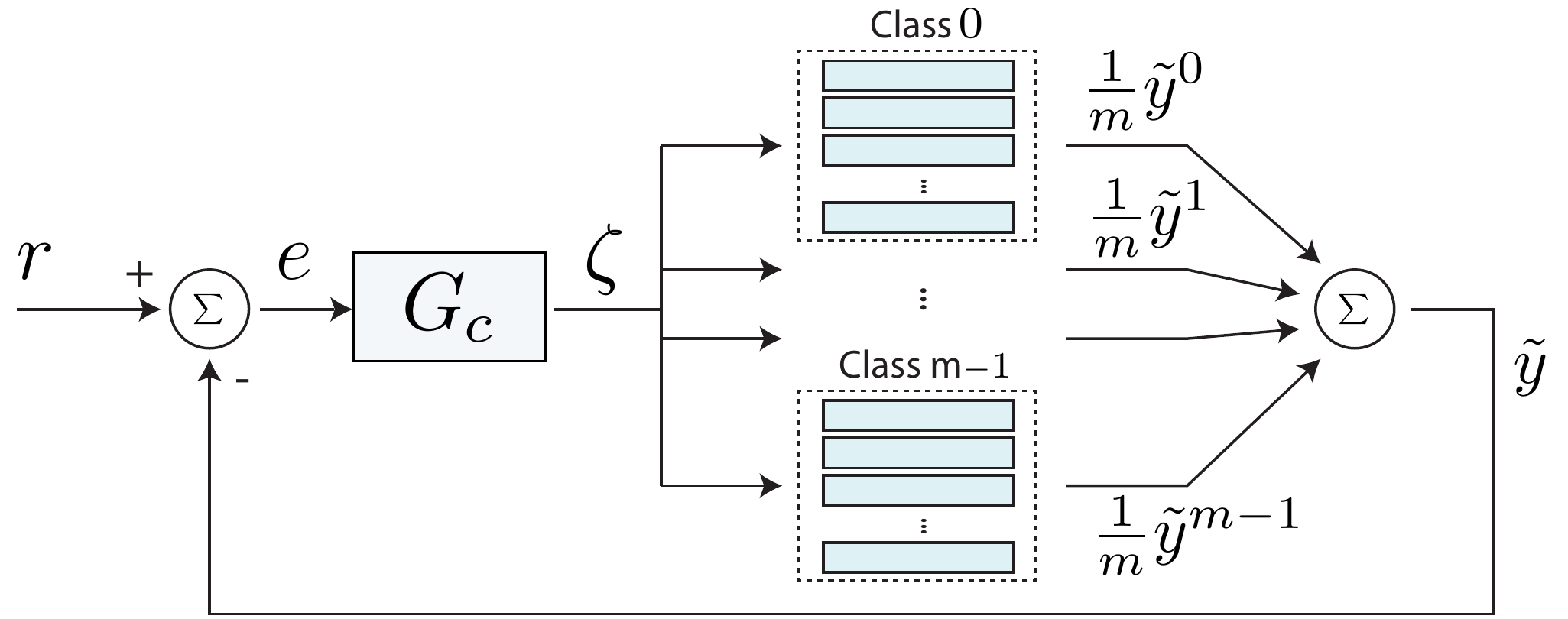}
\vspace{-.2cm}
\caption{The discrete time control signal $\{\zeta_t\}$ is
broadcast at times $t_i = iT_g$, $i=0,1,2,\dots$
in the super-sampled model. A load in class $k$ reads the control signal and makes a decision at consecutive times $\{t_i : i = kT_g +iT,\ i\ge 0\}$. }
\label{f:ss} 
\end{figure} 

To simplify notation the following conventions are made throughout the analysis in this paper:  We henceforth consider a normalized discrete-time model in which $T_g=1$.   The time index $t\in\nat$ refers to time at the grid level,  and $\tau\in\nat$ time for an individual load. 
\spm{NEW:  Yue, let's discuss}

\Fig{f:ss} illustrates  the resulting system architecture.  A linear approximation is again possible, in which the delay is $T_g=T/m$ rather than $T$.

\spm{Gap!  We need to explain that $\tau=mk+i$, and $\tau+1$ is local time, so $\tau+1=mk+i + m$}

\subsection{Intelligent pools}
\label{s:pool}

The paper \cite{meybarbusyueehr14} on which the present work is based considered an application of this system architecture in which the loads were a homogenous collection of pools.  The motivation for considering pools is the inherent flexibility of pool cleaning, and because the total load in a region can be very large.  The maximum load is estimated at 1GW in the state of Florida.   

In \cite{meybarbusyueehr14}  the state space was taken to be the finite set,
\begin{equation}
\state=\{ (\kappa,i) : \kappa \in \{ \oplus,\ominus\} ,\ i\in \{1, 2, 3, \dots\} \}. 
\label{e:poolstate}
\end{equation}
If $X_{\tau} = (\ominus,i)$, this indicates that the load was turned off and has remained off for $i$ time units, and $X_{\tau} = (\oplus,i)$ represents the alternative that the load has been operating continuously for exactly $i$ time units. 
The utility function is equal to the indicator function that a pool is operating, $\util(x) = \sum_i \ind\{ x = (\oplus, i) \}$.

\begin{figure}[h]
\vspace{.2cm}
\Ebox{1}{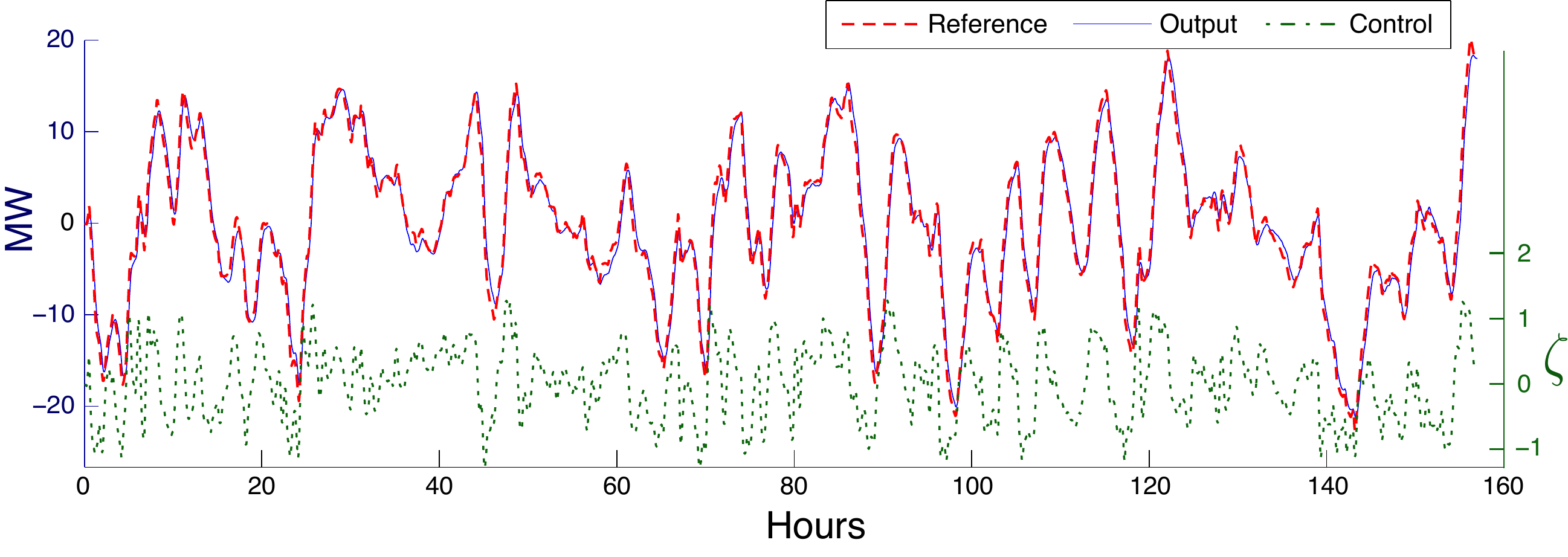}
\vspace{-.05cm}
\caption{The power deviation $\bfmtily$ of a collection of  ``intelligent pools'' tracks nearly perfectly a grid-level power-deviation command $\bfmr$.}
\label{fig:1week_output} 
\vspace{-.25cm}
\end{figure}

\Fig{fig:1week_output} shows results from a simulation experiment similar to those surveyed  in the prior work \cite{meybarbusyueehr14}.   This is based on a symmetric model in which a 12 hour/day cleaning period is desired for each load,  and each load consumes  1kW  when operating.  With $10^5$ pools engaged in providing ancillary service to the grid, tracking at the grid level is  nearly perfect.  A linear control law of the form \eqref{e:zetaGc}
was designed based on the mean field model.  The behavior of the controlled system is as predicted by the mean-field model.

 What was left out in this prior work is any detailed consideration of the quality of service to individual loads.  The pool pump example is ideal for illustrating the possibility of risk for an individual load, and illustrating how this risk can be reduced or even eliminated. 
 

Consider the following measure of QoS for an individual load.  
For a time-horizon $T_f$,  and  function $\health\colon\state\to\Re$,
\begin{equation}
\Health_\tau = \sum_{i=0}^{T_f}  \health(X_{\tau-i}).
\label{e:HealthFinite0}
\end{equation}

In application to the pool pump model we consider the total operation time, so that
\begin{equation}
\health(X_\tau) = \tau_s \sum_j \ind\{X_{\tau}= (\oplus, j) \}
\label{e:health0}
\end{equation}
where $\tau_s = 0.5$ hour is the sampling interval for each pool.

\begin{figure}[h]
\Ebox{.9}{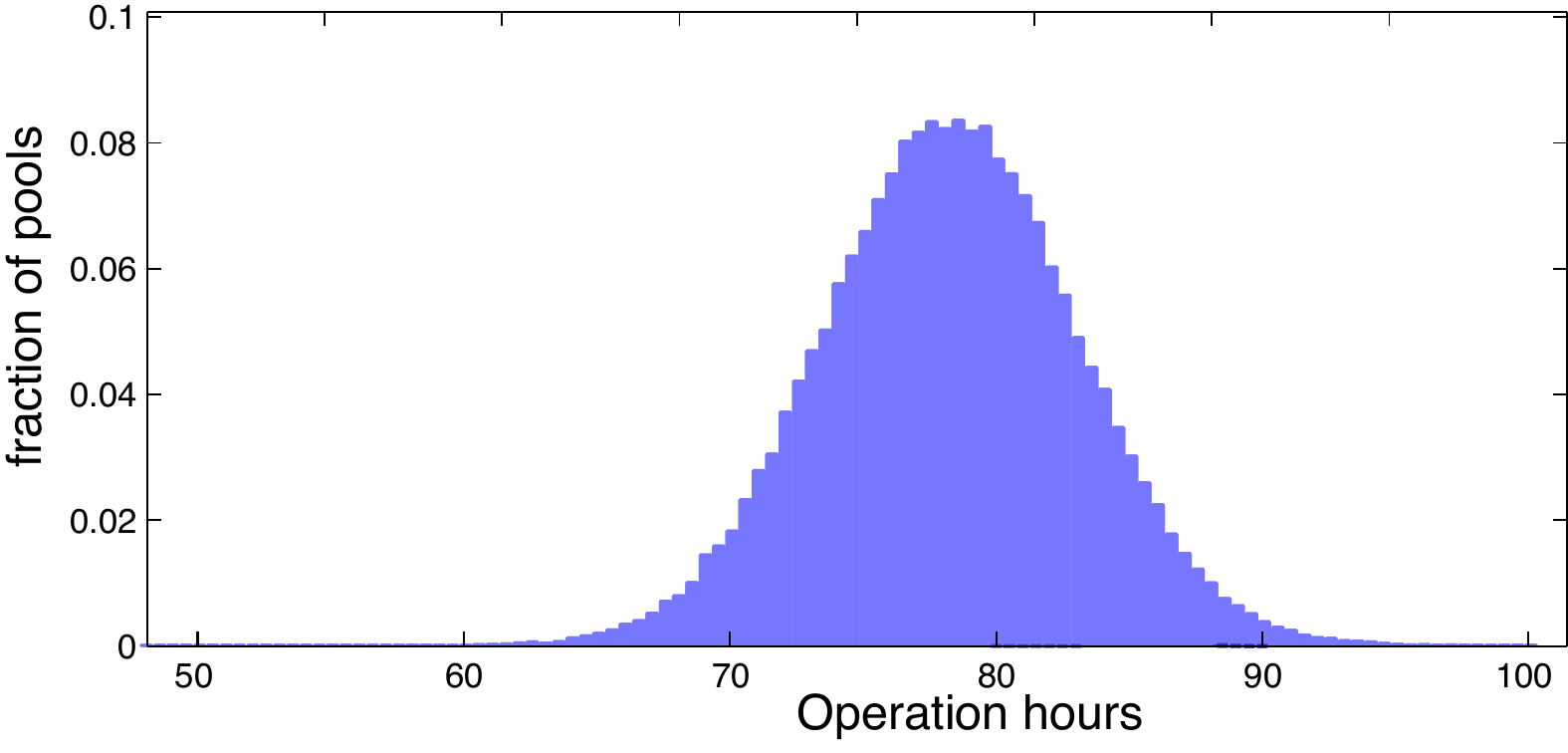}
\vspace{-.2cm}
\caption{Histogram of the moving-window QoS metric \eqref{e:HealthFinite0}.}
\label{fig:1week_hist}
\end{figure}

In the same experiment conducted to produce  \Fig{fig:1week_output},  a histogram of the QoS function \eqref{e:HealthFinite0} was constructed based on the $10^5$ pools.
The time horizon $T_f = 314 $ was chosen corresponding to the reference signal in \Fig{fig:1week_output}, which shows
157 hours of data.  The histogram shown in \Fig{fig:1week_hist} is based on the outcomes for the simulated pools.  
 
The histogram appears Gaussian, with a mean of approximately  77.5 hours ---  this is consistent with the time horizon used in this experiment, given the nominal 12 hour/day cleaning period.  It is evident that a substantial fraction of pools are over-cleaned or under cleaned by 20 hours or more. 
 
An analysis of QoS is presented in the next section based on a model of an individual load in the mean field limit.

\textit{The proofs of the main results   can be found in the Appendix.}

\section{Mean field model for an individual load}
\label{s:MFMload}

In the mean-field limit, the aggregate dynamics are deterministic, following the discrete-time nonlinear control model \eqref{e:MFM}. The behavior of each load remains probabilistic. 

We define the mean field model for \textit{one load} by replacing $(\mu^N_t,\zeta^N_t)$ with its mean-field limit $(\mu_t,\zeta_t)$. 
The justification is that we have a very large number of loads, but our interest is in the statistics of an individual. 

Under the conditions under which the mean field model is obtained as a limit,  the signal $\bfzeta$  is a function of the regulation signal $\bfmr$ that can be expressed as causal feedback,
\[
\zeta_t = \phi_t(\mu_0,\dots,\mu_t, r_0,\dots, r_t),\qquad t\ge 0,
\]
where in the mean field model, the sequence of probability measures evolve according to \eqref{e:MFM}.
Hence $\zeta_t $ is a nonlinear function of $r_0,\dots, r_t$, and the initial condition $\mu_0$ 
(which is assumed to be deterministic).

In some of our analysis it is assumed that $\bfzeta$ is a stationary stochastic process.  The construction of a stationary model for this input process is based on error feedback of the form \eqref{e:zetaGc}, where $y_t=\sum_x \mu_t(x)\util(x)$,  and $\tily_t$ is the deviation from nominal.   To approximate the statistics of $\bfzeta$ we use the linear state space model \eqref{e:LSSmfg} whose output $\bfgamma$ is intended to approximate $\bfmtily$.
The transfer function from  $\bfzeta$ to $\bfgamma$ is denoted $G_p$.  Based on this linear approximation,  we obtain a linear approximation  of the mapping $\bfmr \to \bfzeta$ via the transfer function 
 $G_c/(1+G_cG_p)$    (see any undergraduate textbook on classical control).    

As in our prior work \cite{meybarbusyueehr14,haolinkowbarmey14}, it is assumed that $\bfmr$ is obtained by filtering the total  regulation signal $\bfmr^0$.   In our approximations, the latter is assumed to be derived from filtered white noise with transfer function $G_{wr}$.  


\begin{figure}
\vspace{.1cm}
\Ebox{.75}{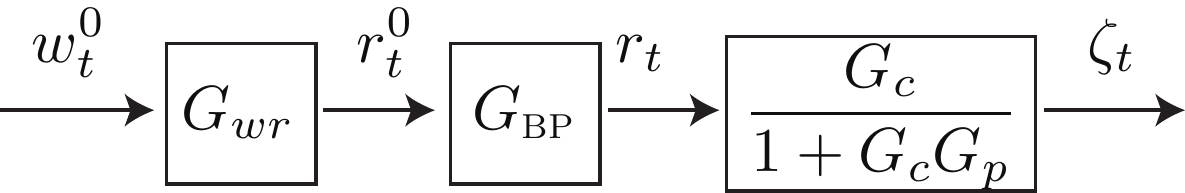}
\vspace{-.2cm}
\caption{The input $\bfzeta$ modeled as a stationary stochastic process}
\label{fig:zetamodel}
\end{figure}  


The super-script $i$ will be dropped in our analysis of a single load. Hence $X(t)$ denotes the non-homogeneous Markov chain whose transition probabilities are defined consistently with \eqref{e:Pzeta},
\begin{equation}
\Prob\{X_{\tau+1} = x' \mid X_{\tau} = x,\ r_0,\dots,r_\tau \} = P_{\zeta_\tau}(x,x'). 
\label{e:Pzetab}
\end{equation}
The construction of the mean field model \eqref{e:MFM} presented in \cite{meybarbusyueehr14} is based on lifting the state space from the $d$-element set
$\state = \{x^1,\cdots,x^d\}$, to the $d$-dimensional simplex $\Spx$. For the $i^{th}$ load at time $\tau$, the element $\piload_{\tau} \in \Spx$ is the degenerate distribution whose mass is concentrated at $x$ if $X_\tau= x$;
that is, $\piload_{\tau} =\delta_x$.

With this state description, the load evolves according to a random linear system, similar to the deterministic dynamics \eqref{e:MFM}:
\begin{equation}
	\piload_{\tau+1}=\piload_{\tau}G_{\tau+1} 
\label{e:piG}
\end{equation}
in which $\piload_{\tau}$ is interpreted as a $d$-dimensional row vector.   The  $d\times d$ matrix 
$G_{\tau}$ has entries $0$ or $1$ only, with $\sum_{x'\in\state} G_{\tau}(x,x')=1$ for all $x\in\state$. It is conditionally independent of $\{\piload_0,\cdots,\piload_{\tau}\}$, given $\zeta_\tau$, with
\begin{equation}
\label{e:EG=P} 
 	\Expect[G_{\tau+1}|\piload_0, \cdots, \piload_{\tau}, \zeta_\tau]=P_{\zeta_\tau}.
\end{equation} 

We denote by $\{X_\tau^\bullet,\ \piload_{\tau}^\bullet : -\infty<\tau<\infty\}$  the two stationary processes obtained with $X_0^\bullet \approx\pi_0$ and $\zeta_t\equiv 0$  for each $t$.  
\spm{NEW:  Yue, I believe $t$ is correct}

\subsection{QoS dynamics and opt-out control}
\label{s:optout}

The local control used to maintain constraints on QoS is a simple ``opt-out'' mechanism.  We formulate a QoS metric $\Health_\tau$ similar to  \eqref{e:HealthFinite0}, along with  upper and lower bounds  $b_+$ and $b_-$.   A load ignores a command to switch state if it will result in $\Health_{\tau+1}\not\in [b_-,b_+]$,  and  take an alternative action so that
$\Health_{\tau+1} \in [b_-,b_+]$.  This ensures that the QoS metric of each load remains within the predefined interval for all time.

The formulation of QoS as the sum 
\eqref{e:HealthFinite0}
is not easily analyzed.  It is more convenient to consider the following discounted sum,
\begin{equation}
\Health_\tau = \sum_{k=0}^\tau \beta^k \health(X_{\tau-k})
\label{e:Health}
\end{equation}
where  the discount factor satisfies $\beta\in(0,1)$, with $\beta\approx 1$.  In simulations presented in \Section{s:num}, the discount factor is chosen so that $\beta^n\approx 1/2$ when $n$ corresponds to a one-week time-horizon.
 
The main goal of this section is to estimate the statistics of $\{\Health_\tau\}$ in steady-state.  For $\beta\approx 1$, a Gaussian approximation is appropriate, so that it is sufficient to estimate its power spectral density  $\psd_{\Health}$.

\begin{figure}
\vspace{.1cm}
\Ebox{.95}{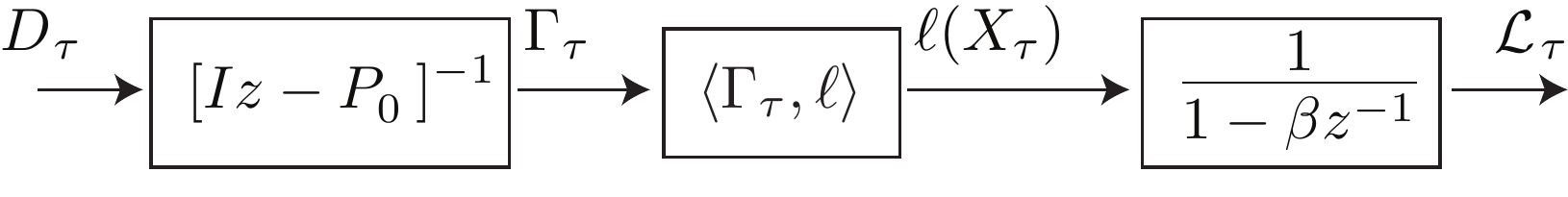}
\vspace{-.2cm}
\caption{Linear model of QoS.} 
\label{f:health}
\end{figure}

The power spectral density   is estimated through several steps:
\begin{romannum}
\item 
We consider the linearized dynamics, so that \eqref{e:piG} is approximated by the linear state space model,
\spm{NEW:  In journal version we must be more precise}
\begin{equation}
\piload_{\tau+1}= \piload_\tau P_0 +  D_{\tau+1}.
\label{e:Sys_delta_lin}
\end{equation}

\item
The disturbance $\bfmD$ in \eqref{e:Sys_delta_lin} is
						approximated by
 the sum of two uncorrelated  components, whose power spectral densities are approximated in \Proposition{t:pi-dynamics} and \Proposition{t:pi-dynamicsSBzeta}.

\item
The approximation in (ii) is based on  a spectral analysis of the input $\bfzeta$, combined with an analysis of the Markov model obtained with $\bfzeta\equiv 0$.

\item
Lastly, given a trace of data $\bfmr$ we obtain an approximation of its power spectral density, and from this we obtain an estimate of the power spectral density of $\bfzeta$.
\end{romannum}
Accepting these approximations, the power spectral density for $\{\Health_\tau\}$ can be obtained via the schematic shown in \Fig{f:health}.



Based on the linear model
\Fig{f:health}
and an approximation for the power spectral density of $\bfmD$, we arrive at an estimate of the power spectral density 
$\psd_{\Health}$.   This admits a stable factorization,
\[
 \psd_{\Health}(\theta) = G_\Health(e^{j\theta}) G_\Health(e^{-j\theta})^\transpose,
\]
where $G_\Health$ is a stable spectral factor, which is taken to be a row-vector of transfer functions. On letting $\{g_k\}$ denote its impulse response,
$
 G(z) = \sum_{k=0}^\infty g_k z^k
$,
the variance of QoS is obtained in two forms,
 \[
\sigma^2_{\Health} = \frac{1}{2\pi} \int_0^{2\pi} \psd_{\Health}(\theta) \, d\theta = \sum_{k=0}^\infty \|g_k\|^2.
\]

We now proceed with the construction of this linear model.

\subsection{LTI approximations}

The random linear system \eqref{e:piG} can be described as a linear system driven by ``white noise'':
\begin{equation}
	\piload_{\tau+1}=\piload_{\tau}P_{\zeta_\tau}+\Delta_{\tau+1}
	\label{e:Sys_delta}
\end{equation}
where, $\{\Delta_{\tau+1}=\piload_{\tau}(G_{\tau+1} -P_{\zeta_\tau}): \tau\geq 0 \}$ is a martingale difference sequence.
We can compute the covariance of the noise in one special case: 
\begin{proposition}
\label{t:pi-dynamics}
Consider the stationary sequence $\{\Delta_{\tau}^\bullet\}$ obtained for the system in which  $\piload_{\tau} = \piload_{\tau}^\bullet$ for all $\tau$.
Its covariance   is given by,
\begin{equation}
\Sigma^\Delta \eqdef \Cov ( \piload_{\tau}^\bullet[G_{\tau+1}-P_0] ) = \Pi - P_0^\transpose \Pi P_0
\label{e:sigma-delta}
\end{equation} 
where $\Pi=\diag(\pi_0)$.
\qed
\end{proposition}

The LTI approximation of the nonlinear system \eqref{e:Sys_delta} is more complex than 
the linearized mean field model \eqref{e:LSSmfg} because of the presence of randomness.
We have the small-signal approximation as before: applying \eqref{e:Pder},  the approximation   \eqref{e:Sys_delta_lin}  holds, in which the disturbance consists of two terms,
\begin{equation}
D_{\tau+1} \eqdef
B_{\tau}^\transpose \zeta_\tau +\Delta_{\tau+1}
\label{e:D}
\end{equation}
where
  $B_{\tau}^\transpose =\piload_{\tau} \clE$.  The recursion \eqref{e:Sys_delta_lin}  is of the form of the state equation that defines $\bfPhi$ in \eqref{e:LSSmfg}, with the addition of two  ``disturbances''. We can extend \Proposition{t:pi-dynamics}
to approximate  the power spectral density of the overall disturbance, but not in closed form:

\begin{lemma}
\label{t:pi-dynamics2}
The power spectral density of the stationary sequence $\{\zeta_\tau  \piload_{\tau}^\bullet \clE +\Delta_{\tau+1}^\bullet \}$ can be expressed as the sum,
\[
\psd(\theta) = \psd_{B\zeta}(\theta) + \psd_\Delta(\theta),
\]
where $ \psd_\Delta(\theta) = \Sigma^\Delta$ for all $\theta$. The power spectral density $ \psd_{B\zeta}(\theta)$ for $\piload_{\tau}^\bullet \clE \zeta_\tau$ is the Fourier transform of the product of the autocorrelation functions for $\piload_{\tau}^\bullet \clE$ and $\zeta_\tau$. 
\qed
\end{lemma}

\proof
The random vectors $\zeta_\tau\piload_{\tau}^\bullet \clE$ and $\Delta_{\tau+1}$ are uncorrelated by construction, which is why we obtain a sum. The sequence $\{\Delta_{\tau+1}\}$ is uncorrelated, so that $ \psd_\Delta(\theta) $ is independent of $\theta$. Finally, 
$ B_{\tau}^\bullet \eqdef [ \piload_{\tau}^\bullet \clE]^\transpose$ 
and $ \zeta_\tau $ are independent, which is why we obtain a product of autocorrelation functions. 
\qed

Computation of $\psd_{B\zeta}$ is illustrated in an example.  First we require some additional notation.
The functions $\{\varepsilon_k\}$ are defined so that  the $k$th component of $B_{\tau}$ is given by $\varepsilon_k(X_\tau)$:
\[
\varepsilon_k(x^j) =  \clE_{j,k},\quad 1\le j,k\le d.
\]
For any $\poleshift\in\Co$ satisfying $|\poleshift|<1$, the resolvent matrix is the discounted sum,
\[
U_\poleshift= \sum_{n=0}^\infty \poleshift^n P^n.
\]
 For two functions $f,g\colon\state\to\Re$ denote,
\[
\langle f,g\rangle = \sum \pi_0(x) f(x) g(x).
\]

\begin{proposition}
\label{t:pi-dynamicsSBzeta}
Consider the special case in which the autocorrelation function for $\bfzeta$   is given by
\begin{equation}
R_\zeta(n) = \sigma^2_\zeta \rho^{n}, \quad n\ge 0,
\label{e:Rzeta}
\end{equation}
where  $\rho\in(0,1)$ and $\sigma^2_\zeta>0$. The $k,l$ entry of the power spectral density matrix is then given by,
\begin{equation} 
\psd_{B\zeta}^{kl} (\theta) =\bigl(\langle \varepsilon_k  ,U_\poleshift \varepsilon_l \rangle  +\langle\varepsilon_l , U_\poleshift \varepsilon_k \rangle ^* - \langle \varepsilon_k,  \varepsilon_l\rangle   \bigr) \sigma^2_\zeta 
%
\label{e:psdBzeta}
\end{equation}
where $\poleshift =\rho^m e^{-j\theta }$, and the superscript ``$*$'' denotes complex conjugate.
\qed
\end{proposition}

\spm{There was a gap here:   We were forgetting that we are sampling $\bfzeta$. Hence we need the autocorrelation of the subsampled process $\{\zeta_{km}:  k\in\nat\}$.  This is why I now have   $\poleshift =\rho^m e^{-j\theta }$}

  It is not difficult to extend the result to the case where $\bfzeta$ has a power spectral density, with stable spectral factor $G_\zeta$:
\[
 \psd_{\zeta}(\theta) = G_\zeta(e^{j\theta}) G_\zeta(e^{-j\theta}) .
\]
If the transfer function $G_\zeta(z)$ has distinct poles $\{\rho_1,\dots,\rho_{n_G}\}$,  then we can construct constants
$\{A_1,\dots,A_{n_G}\}$ such that  for $n\ge 0$,
\[
R_\zeta(n) = \sum_{j=1}^{n_G} A_j \rho_j^n.
\]
Based on this representation, the power spectral density $\psd_{B\zeta}$ can be expressed as a sum of $n_G$ terms, each of the form given in 
\eqref{e:psdBzeta}.

\subsection{Model for the regulation signal.}

This section is concluded with discussion on the construction of the transfer function $G_{wr}$ shown in \Fig{fig:zetamodel}.

\begin{figure}
\vspace{.2cm}
\Ebox{1}{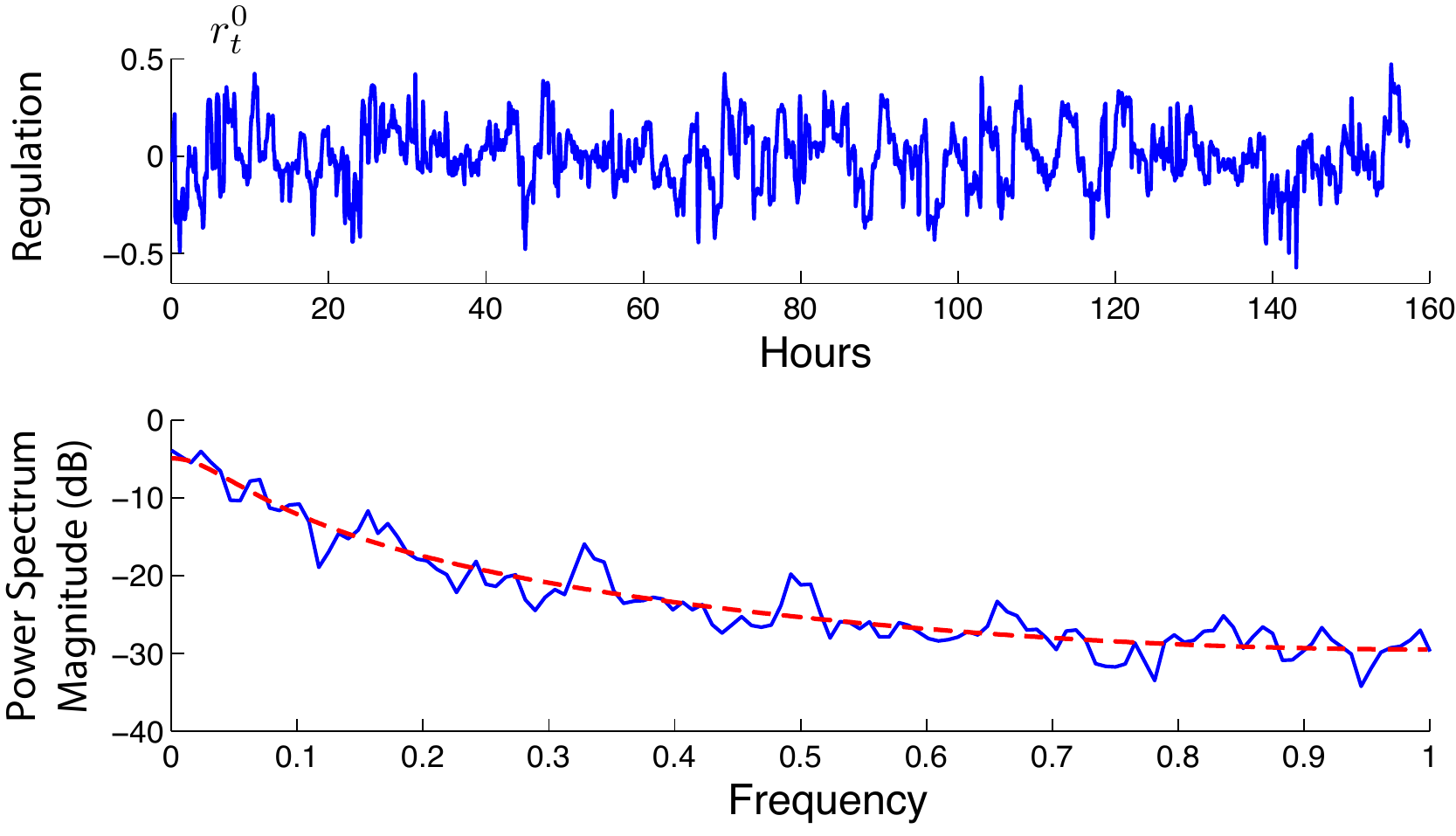}
\vspace{-.2cm}
\caption{Scaled BPA signal.}
\label{fig:BPA}
\vspace{-.2cm}
\end{figure}

\Fig{fig:BPA} shows the normalized regulation signal for a typical week within the Bonneville Power Authority (BPA)~\cite{BPA}.   The power spectral density plot shown was estimated using Matlab's  {\tt psd} command.  


To obtain a statistical model,  it is assumed that $\bfmr^0$ evolves as the ARMA model,
\begin{equation}
	r_t^0 +a_1r_{t-1}^0+a_2r_{t-2}^0=w_t^0+b_1w_{t-1}^0
 \label{e:ARMA}
\end{equation}
in which $\bfmw^0$ is white noise with variance $\sigma_w^2$. 
The extended least squares (ELS) algorithm was used to estimate the coefficients $a_1$, $a_2$, $b_1$, and the variance $\sigma_w^2$. Based on data samples in \Fig{fig:BPA}, the ELS algorithm terminated at $[a_1, a_2, b_1]^T=[-0.9009, 0.0365, 0.0859]^T$, and $\sigma_w^2=0.005$. In the $z$-domain, its transfer function is expressed as below,
\begin{equation}
	G_{wr}(z)=\frac{1+0.08594z^{-1}}{1-0.9009z^{-1}+0.03653z^{-2}}.
 \label{e:G_{BPA}(z)}
\end{equation}
The dashed line in \Fig{fig:BPA} is the estimate of the spectrum given by $ |R^0(e^{jw})|^2 = \sigma_w^2 |G_{wr}(e^{jw})|^2$.



\section{Application to intelligent pools}
\label{s:num}

Numerical experiments were conducted on the pool pump model described in \Section{s:pool} to investigate the performance of local opt-out control. We  show through simulations the improvement on QoS and the impact on the grid with this additional local control. In addition, we provide an analytic formula which can be used to obtain real-time estimate of QoS.

\subsection{Simulation setup}

In order to monitor the health of pools over a long time horizon, a simulated regulation signal was generated using the filter \eqref{e:G_{BPA}(z)}. Driven by white noise $N(0, \sigma_{w}^2)$ with variance $\sigma_{w}^2 = 0.005$ and sampling time $T_g =5$ minutes, we obtained a $400$-hour long simulated regulation signal  $\bfmr^0$.
This signal was then passed through a low-pass filter to generate the reference signal $\bfmr$ used in the experiments.

The supersampling approach was used with $m=6$,   corresponding to a local sampling time equal to $T= mT_g=30$ minutes.  The simulation used $10^5$ homogeneous Markov models;   
each pool pump is assumed to consume $1$~kW during operation.

In the QoS function a normalized 
version of the function \eqref{e:health0}  is used, defined by,
\begin{equation}
\health(X_\tau) = \sum_j [\ind\{X_{\tau}= (\oplus, j) \} - \ind\{X_{\tau} = (\ominus, j) \} ] .
\label{e:health}
\end{equation}
If $\Health_\tau>0$ then the pool has received too much cleaning, and
if $\Health_\tau<0$ then the pool has been under-cleaned.
The acceptable  interval for opt-out control was set to $[-20, \; +20]$ ($b_-=-20$ and $b_+=20$, following the notation at the top of \Section{s:optout}).

The discount factor  in the QoS function \eqref{e:Health} was chosen as $\beta = 0.9975$.

%

 \spm{let's eliminate this subsection title in revision}


\subsection{Performance of opt-out control}

The QoS without opt-out control is illustrated in the histogram shown on the left hand side of \Fig{fig:BothHist}, based on samples of the  QoS function \eqref{e:Health} from each load over $400$ hours. A Gaussian distribution with mean zero and variance $\sigma_{\Health}^2 $ is plotted as the dashed line in this figure,  based on the variance approximation obtained in \Section{s:MFMload}.  The approximation is close to the simulation result.

\begin{figure} 
\vspace{.2cm}
\Ebox{.95}{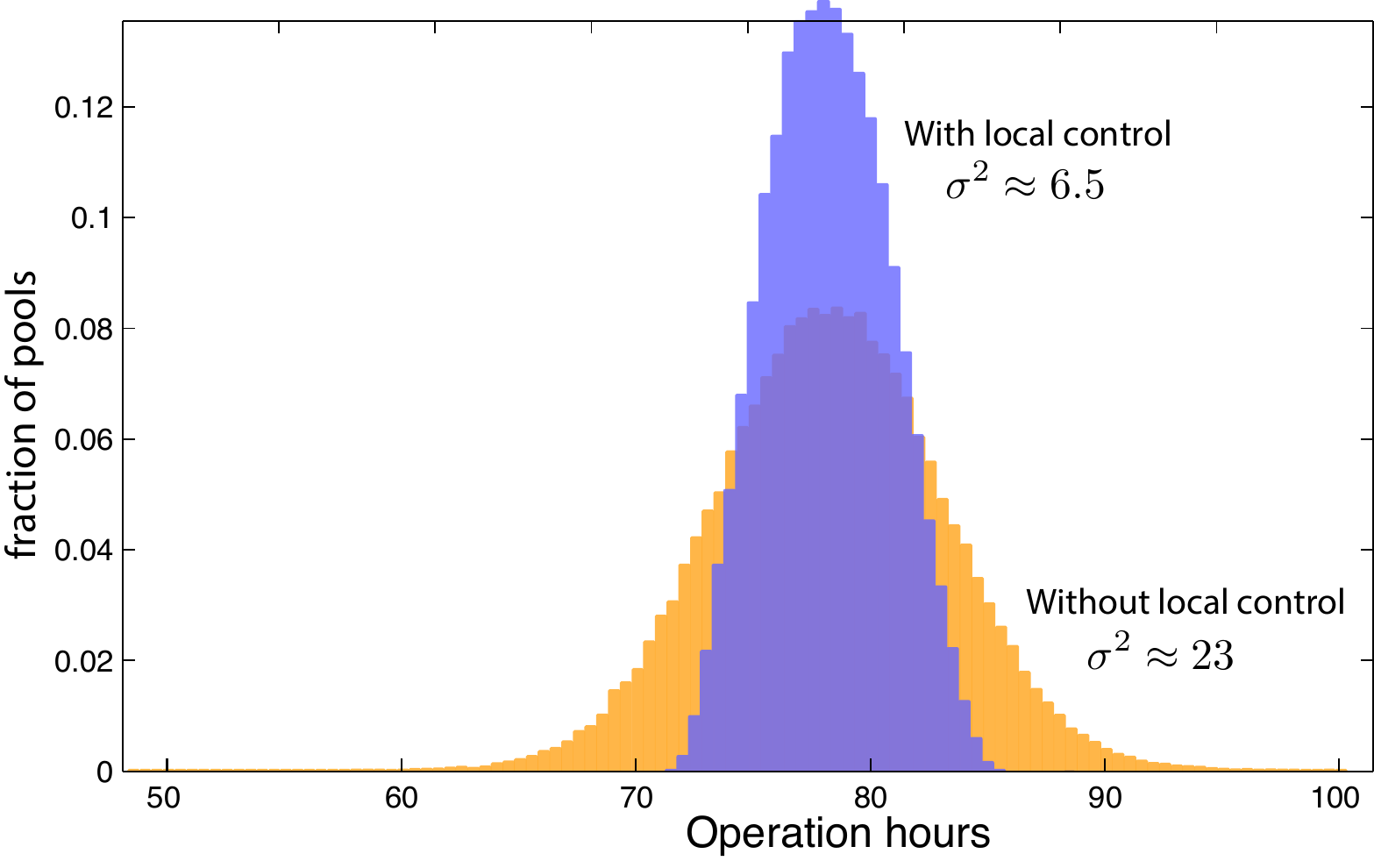}  
\vspace{-.2cm}
\caption{Comparison of the histograms for the \textit{moving-window} QoS indicator \eqref{e:HealthFinite0}, with and without local control.  }
\label{fig:OneWeekHistBoth} 
\end{figure}

The regulation signal used to generate  on the left hand side of \Fig{fig:BothHist} was re-used in experiments to investigate the impact of opt-out control on the grid, and on the QoS of an individual load. 

QoS is dramatically improved with opt-out control:
The histogram shown on the right hand side of \Fig{fig:BothHist} is truncated to the interval $[-20, \;+20]$ as desired.  The dashed line is the conditional Gaussian density, truncated to this interval, with the same mean and variance as the Gaussian distribution  shown on the left hand side of this figure.  The conditional density is a good fit to the histogram obtained through simulations.

These plots are all based on the discounted QoS metric \eqref{e:Health}.  
\Fig{fig:OneWeekHistBoth} shows a comparison of two histograms of the moving-window QoS metric \eqref{e:HealthFinite0}, with and without local control.   The histogram with larger variance is exactly the same as shown in 
\Fig{fig:1week_hist}.  The variance is reduced by more than a factor of three using local control.

\begin{figure}
\vspace{.2cm}
\Ebox{1}{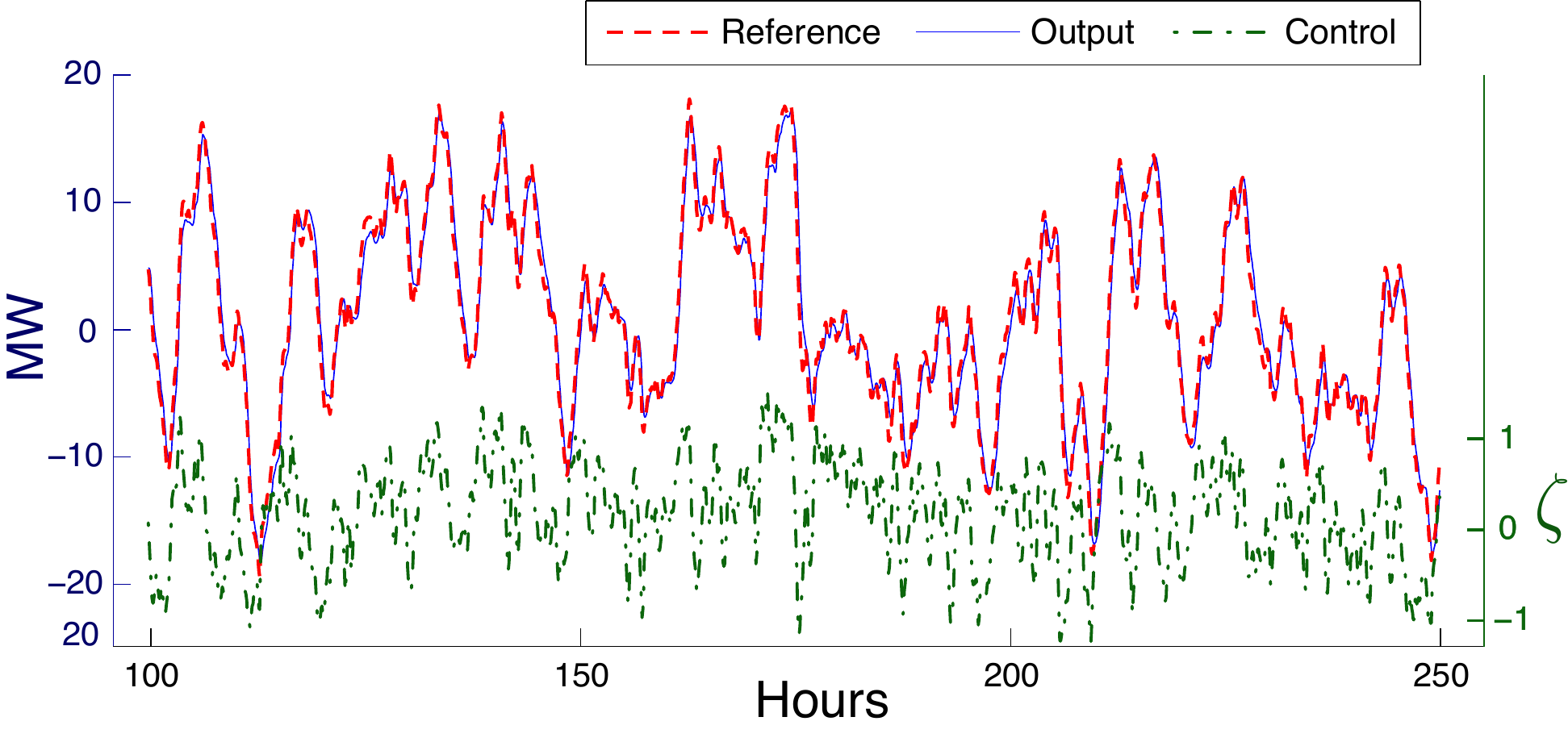}
\vspace{-.2cm}
\caption{Tracking performance with local control.}
\label{fig:LC_r1_output} 
\vspace{-.2cm}
\end{figure}

The grid level tracking performance shown in \Fig{fig:LC_r1_output} remains nearly perfect. This is surprising, given the  improvement  in QoS. The explanation is that very few loads opt out:  In \Fig{fig:LC_r1_optout} we see that in this simulation, no more than $3\%$ of loads opt out at any time.


There are limitations on the capacity of ancillary service from  a collection of loads, and it is likely that op-out control will reduce capacity. We next obtain an approximation of QoS to better understand the grid-level impact of opt-out control.

\spm{Figures are not normalized}

\subsection{QoS analysis}

The variance of QoS was approximated  in \Section{s:MFMload} based on a small-signal analysis.

The following proposition gives an approximation of its mean at any time.  This requires bridging time at the grid and load levels:  Any time $t$ can be expressed as $t = m\tau+j$, where $m$ is the   super-sampling parameter, $\tau \in \nat$, and $j \in \{0,1, \cdots, m-1\}$.   The average  QoS at time $t$ is denoted, 
\[
\bar{\Health_t} = \frac{1}{N}\sum_{i=1}^{N}\Health_t^i
\]  
where $N$ is the number of pools, $\Health_t^i$ is the QoS for load $i$ at time $t$,  defined using grid-level time.   For each $i$, the function $\{\Health_t^i : t=0,1,...\}$ is piecewise constant, with potential jumps spaced $T_g$ minutes apart. 

A   discounted reference sequence is defined by,
\begin{equation}
R_t^{\beta} = \sum_{k=0}^{\tau} \beta^k r_{m(\tau-k)+j} 
\label{e:R^beta}
\end{equation}  

\begin{proposition}
\label{t:H_approx}
Suppose there is perfect tracking: $\tily_t=r_t$ for all $t$, where $\tily_t=y_t-\bary^0$, and $\bary^0=\frac{1}{2}$. Then,
\begin{equation}
\bar{\Health_t} = 2 R_t^{\beta}
\label{e:Health=2R}.
\end{equation}
\qed
\end{proposition}

\begin{figure}
\vspace{.2cm}
\Ebox{1}{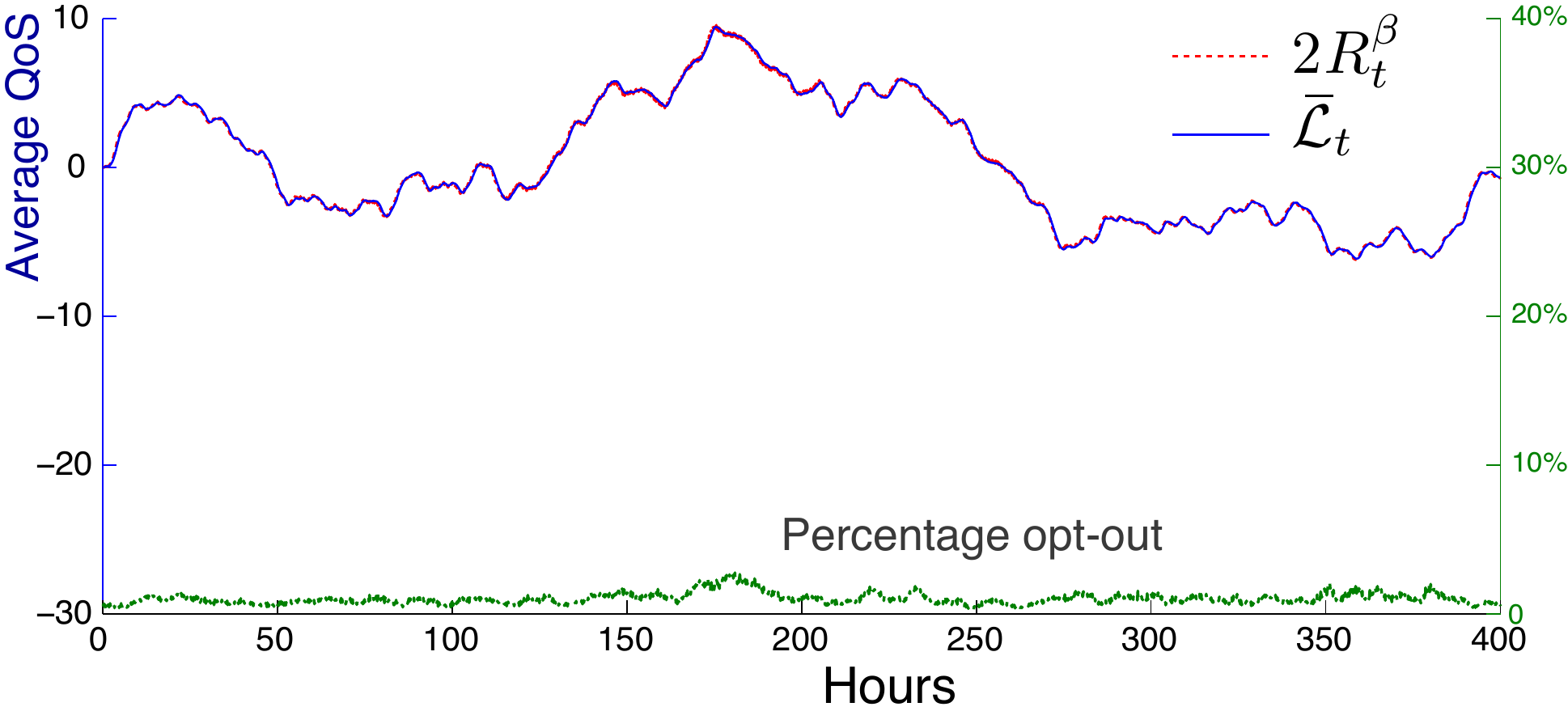}
\caption{Percentage of pools that ``opt out'' as a function of time.}
\vspace{-.2cm}
\label{fig:LC_r1_optout}
\end{figure} 

In practice $\tily_t$ only approximates $r_t$, so \eqref{e:Health=2R} is only an approximation. This is illustrated in \Fig{fig:LC_r1_optout}, where the dotted line is $2R_t^{\beta}$, and the solid line is $\bar\Health_t$. They nearly coincide because of the nearly perfect tracking observed in \Fig{fig:LC_r1_output}.

If the reference signal $\bfmr$  takes a long and large positive or negative excursion, the proposition suggests that many loads will receive poor quality of service, and opt out of service to the grid.

\begin{figure}
\vspace{.2cm}
\Ebox{1}{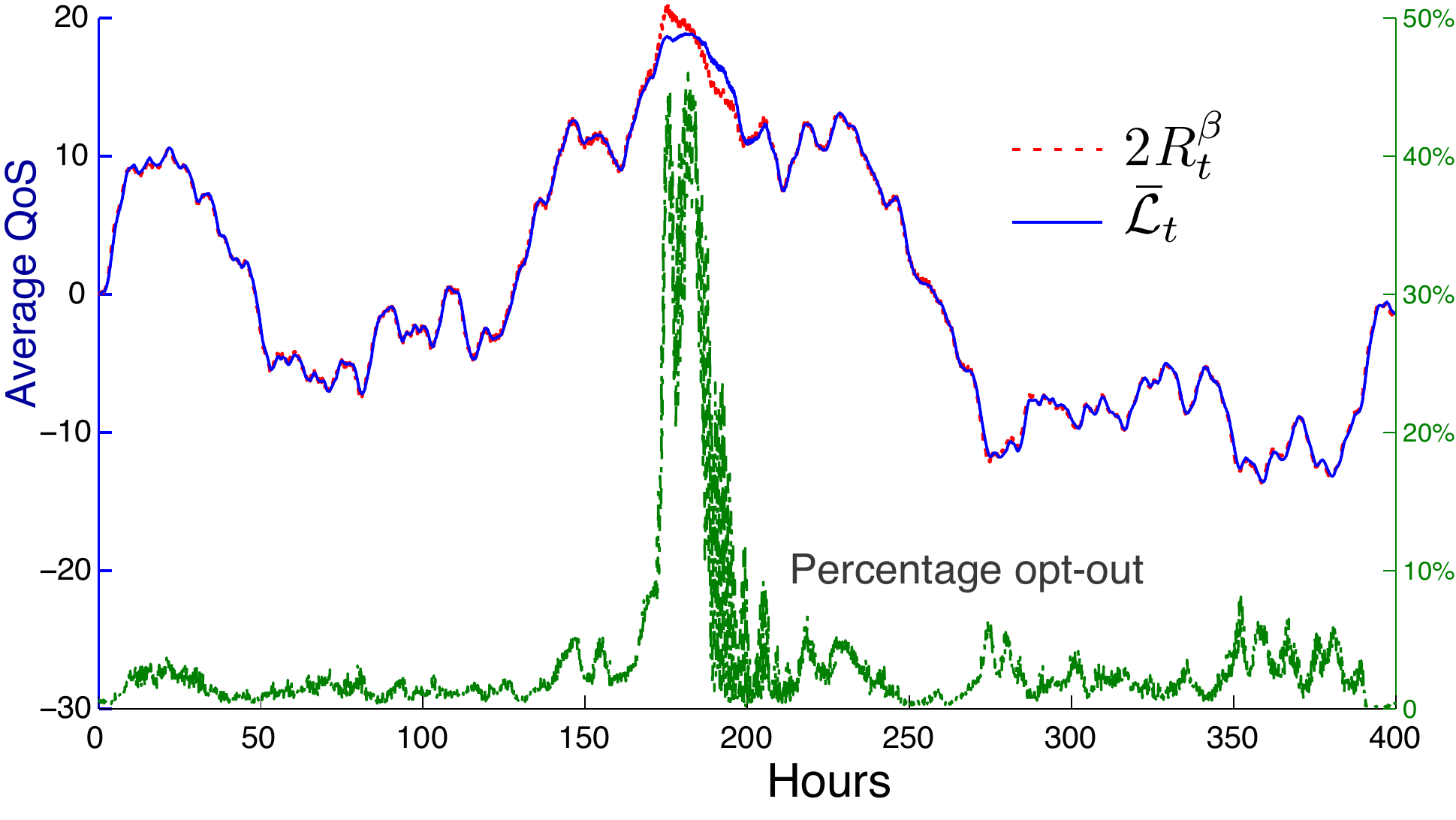}
\vspace{-.2cm}
\caption{Percentage of pools that opt out for a larger reference signal.}
\label{fig:LC_r2_optout}
\end{figure}

\begin{figure}
\Ebox{1}{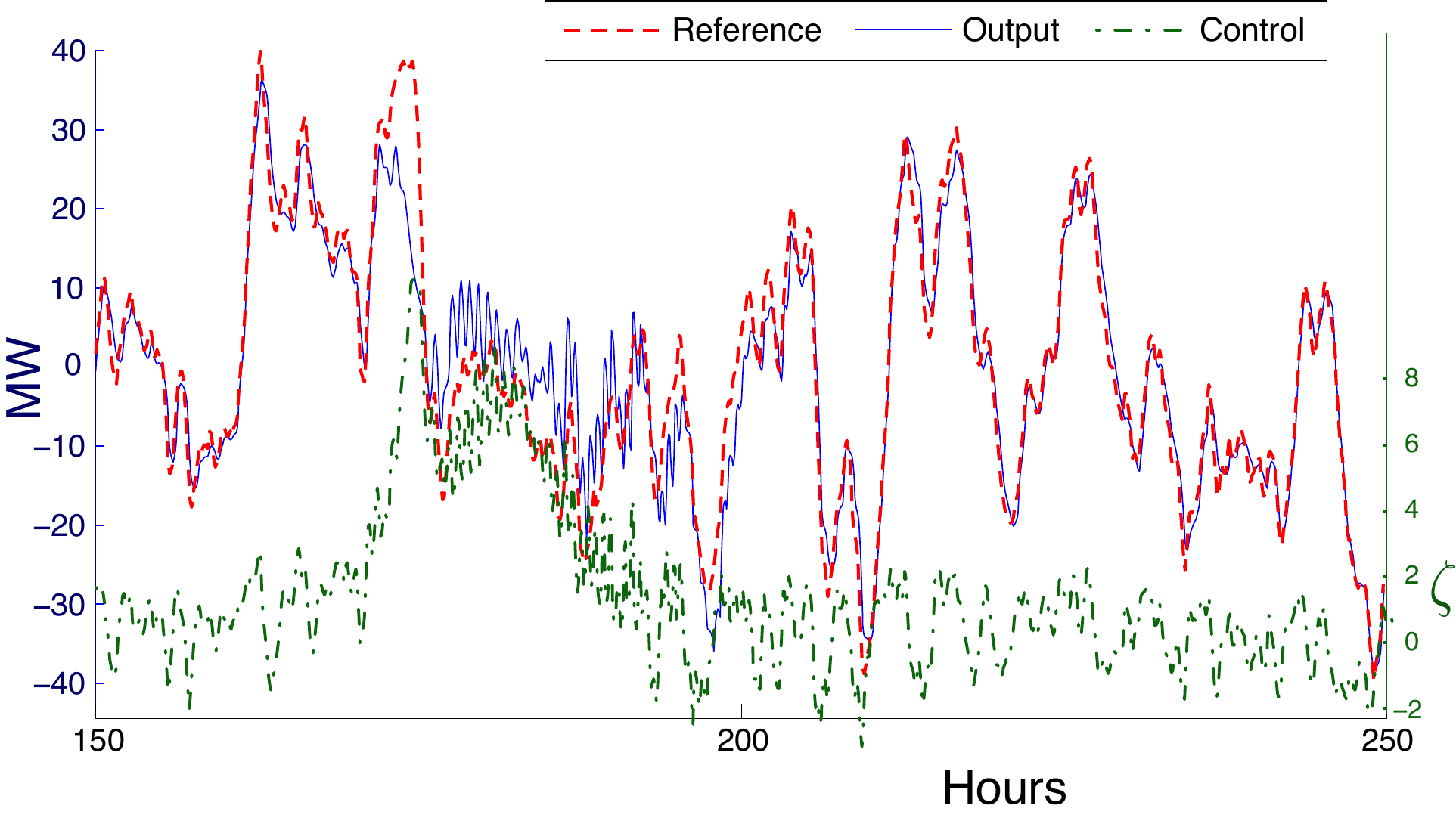} 
\vspace{-.2cm}
\caption{Tracking performance using a larger reference signal with local control --- Close-up over period of poor average QoS seen in \Fig{fig:LC_r2_optout}.}
\label{fig:LC_r2_output} 
\vspace{-.2cm}
\end{figure}

An experiment was conducted to illustrate this conclusion.  The reference signal was scaled by $2.2$, which is approximately the largest scaling possible for the system without opt-out control. Results are summarized in \Fig{fig:LC_r2_optout}. In this situation, $2R_t^\beta$ rose over $20$ at $t\approx 175 hr$, which is the upper bound for an individual load. This caused many loads to opt out, and the corresponding tracking performance degraded at  $t\approx 175 hr$,  as shown in \Fig{fig:LC_r2_output}.

\begin{figure}
\vspace{.3cm}
\Ebox{.9}{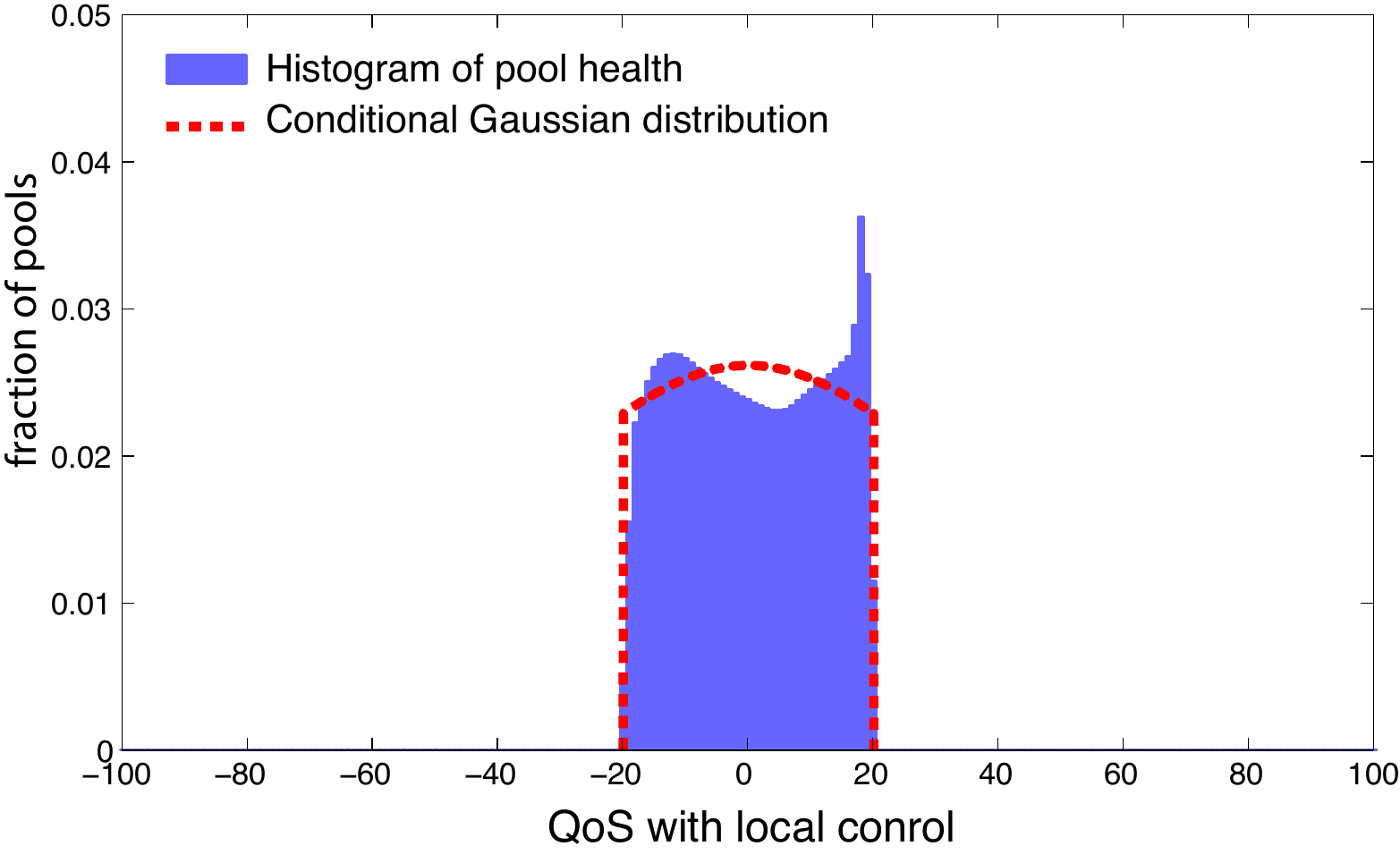} 
\vspace{-.2cm}
\caption{Histogram of QoS subject to a larger reference signal.}
\label{fig:LC_r2_hist}
\end{figure}

The histogram of QoS in \Fig{fig:LC_r2_hist} is still within the interval $[-20, \;+20]$, but it has a greater fraction of pools close to the boundary $+ 20$.


\section{Conclusions}

The main technical contribution of this paper is the approximation of QoS for an individual load.   It is remarkable that it is possible to obtain accurate estimates of first and second order statistics for an individual load, taking  into account second order statistics of exogenous inputs (in this case the reference signal),  along with correlation introduced by the Markovian model.  It is also remarkable that strict bounds on QoS can be guaranteed while retaining nearly perfect grid-level tracking.

The paper has focused on a narrow definition of QoS.  In the case of pools, a consumer will have many constraints besides operation time.  For example, 
the consumer
 may wish to impose limits on the number of mode changes per week.    The basic message of the paper is unchanged:  As long as the system is designed so that opt-out is not too frequent,  the mean-field model will be subject to a small disturbance that can be handled through control design at the grid level.

The linear model used in this work will be useful for other applications, such as online estimation.  We have assumed that the BA observes  the output $\bfmy$,  perhaps with some measurement error.  Based on the linear model we can apply the Kalman filter to obtain an estimate of the distribution of loads,
so that the BA can obtain online approximations of  variables that are linear functions of the empirical distribution~\eqref{e:empDist}.

In particular,   state estimation may be valuable for estimating potential capacity of demand response, which will vary over time. Recall the plot shown in \Fig{fig:LC_r2_output}, illustrating how tracking performance degrades when a large number of loads opt-out of service to the grid. A Kalman filter can provide real-time  estimates of the number of loads that have opted-out, which directly translates to a capacity bound.

\bibliographystyle{IEEEtran}
\bibliography{strings,markov,q}

 \appendix

\subsection{Proof of \Proposition{t:pi-dynamics} }

The random variables $\piload_\tau$ and $G_{\tau+1}$ are independent, with respective means $\pi_0$ and $P_0$.  
The row vector $\piload_\tau$ has entries zero or one, and the same is true for the matrix $G_{\tau+1}$.

The $k,l$ entry of the covariance matrix $\Sigma^\Delta$    can be computed as follows.  First, using the fact 
that $ \piload_\tau(i) \piload_\tau(j) =0$ if $i\neq j$,
\[
\begin{aligned}
&\Sigma^\Delta_{k,l} \\
&=\Expect\Bigl[  \Bigl( \sum_i \piload_\tau(i) [G_{\tau+1}(i,k) - P_0(i,k) ] \Bigr) \\
&\qquad \quad \Bigl( \sum_j\piload_\tau(j) [G_{\tau+1}(j,l) - P_0(j,l) ] \Bigr) \Bigr]
 \\
 &
 =
 \Expect\Bigl[  \sum_i \piload_\tau(i)  [G_{\tau+1}(i,k) - P_0(i,k) ]   [G_{\tau+1}(i,l) - P_0(j,l) ] \Bigr]
  \\
 &
 =\sum_i \pi_0(i)
 \Expect\Bigl[    [G_{\tau+1}(i,k) - P_0(i,k) ]  [G_{\tau+1}(i,l) - P_0(i,l) ] \Bigr]
\end{aligned}
\]
where the last line uses independence of $\piload_\tau$ and $G_{\tau+1}$.  Further manipulations give, 
\[
\begin{aligned}
\Sigma^\Delta_{k,l} 
	& =\sum_i \pi_0(i)   \Bigl( \Expect\bigl[    G_{\tau+1}(i,k)   G_{\tau+1}(i,l) \bigr]  
\\
	& \qquad \qquad \qquad  \qquad \qquad  \qquad        - P_0(i,k)P_0(i,l) \Bigr)
 \\
 	&   =\sum_i \pi_0(i)  \Bigl( \Expect\bigl[       G_{\tau+1}(i,k) \bigr] \ind\{k=l\} 
\\
	& \qquad \qquad \qquad  \qquad \qquad  \qquad          - P_0(i,k)P_0(i,l) \Bigr)
 \\
 	&   =\sum_i \Bigl(\pi_0(i)  P_0(i,k)  \ind\{k=l\}  - P_0(i,k)P_0(i,l) \Bigr) 
\end{aligned}
\]
where we have used the property that $  G_{\tau+1}$ has Bernoulli entries.  The right hand side is precisely  \eqref{e:sigma-delta}.
\qed

\subsection{Proof of \Proposition{t:pi-dynamicsSBzeta}.}
 
We set $\sigma^2_\zeta =1$ to simplify  notation.

The power spectral density matrix $\psd_{B\zeta}(\theta)$ is the Fourier transform of the autocorrelation matrix for $\bfmB \bfzeta$.  
The stochastic processes  $\bfmX=\bfmX^\bullet$ and $ \bfzeta $ are independent under the conditions of the proposition.  Consequently,  the $k,l$ entry of $\psd_{B\zeta}(\theta)$ can be represented as,
\[
\begin{aligned}
	\psd_{B\zeta}^{kl} (\theta) &= \sum_{\tau=-\infty}^{\infty} \Expect \Bigl[ [\varepsilon_k(X_0)  \zeta_0] [\varepsilon_l(X_{\tau}) \zeta_{m\tau}] \Bigr] e^{-j\tau\theta} \\
						&=\sum_{\tau=-\infty}^{\infty}  \Expect [\varepsilon_k(X_0) \varepsilon_l(X_{\tau}) ]  \Expect[\zeta_0 \zeta_{m\tau}] e^{-j\tau\theta} \\ 
						&= \sum_{\tau=-\infty}^{\infty} \Expect [\varepsilon_k(X_0) \varepsilon_l(X_{\tau})  ]   \rho^{|m\tau|} e^{-j\tau\theta} 
\end{aligned}
\]
where we have used $\sigma^2_\zeta=1$, and each expectation is in steady-state.

The right hand side is expressed as the sum of  three terms:
\begin{equation}
\begin{aligned}
	\psd_{B\zeta}^{kl} (\theta) &=   \sum_{\tau=0}^{\infty} (\rho^m e^{-j\theta})^\tau \Expect [\varepsilon_k(X_0) \varepsilon_l(X_{\tau}) ]  \\
						& \quad +   \sum_{\tau=-\infty}^{0} (\rho^m e^{j\theta})^{-\tau} \Expect [\varepsilon_k(X_0) \varepsilon_l(X_{\tau}) ]  \\
						& \quad -   \Expect [\varepsilon_k(X_0) \varepsilon_l(X_0) ] 
\label{e:3terms}
\end{aligned}
\end{equation}
The last term in this expression is the inner product $   \Expect [\varepsilon_k(X_0) \varepsilon_l(X_0) ]  =    \langle \varepsilon_k,  \varepsilon_l\rangle $
that appears in \eqref{e:psdBzeta}.

Using the smoothing property of conditional expectation and the definition $\poleshift =\rho^m e^{-j\theta }$, the first term in \eqref{e:3terms} can be expressed in terms of the resolvent:
\[
\begin{aligned}
     \sum_{\tau=0}^{\infty}  (\rho^m e^{-j\theta})^\tau & \Expect [\varepsilon_k(X_0) \varepsilon_l\,(X_{\tau}) ]   \\
	&=   \sum_{\tau=0}^{\infty} \poleshift^\tau \Expect \Bigl[\varepsilon_k(X_0) \Expect [\varepsilon_l\,(X_{\tau}) |X_0] \Bigr]  \\
	&=   \sum_{\tau=0}^{\infty} \poleshift^\tau \Expect [\varepsilon_k(X_0) P^{\tau} \varepsilon_l\,(X_0)] \\
	&=   \Expect \Bigl[\varepsilon_k(X_0) \sum_{\tau=0}^{\infty} \poleshift^\tau P^\tau \varepsilon_l\,(X_0) \Bigr]  \\
	&=   \langle \varepsilon_k  ,U_\poleshift \varepsilon_l \rangle 
\end{aligned}
\]
The second term  in \eqref{e:3terms} is transformed into a form similar to the first term,
\[
\begin{aligned}
 \sum_{\tau=-\infty}^{0} (\rho^m e^{j\theta})^{-\tau}& \Expect [\varepsilon_k(X_0) \varepsilon_l(X_{\tau}) ] \\
	&=  \sum_{\tau=0}^{\infty} (\poleshift^*)^{\tau} \Expect [\varepsilon_l(X_0) \varepsilon_k(X_{\tau}) ]  
\end{aligned}
\]
which is equal to $
	  \langle\varepsilon_l , U_\poleshift \varepsilon_k \rangle ^* $.   These representations and \eqref{e:3terms} establish the desired conclusion \eqref{e:psdBzeta}.
\qed

\subsection{Proof of \Proposition{t:H_approx}.}

Denote for any time $t = m\tau+j$, 
\[
		\bar{\health_t} = \frac{1}{N} \sum_{i=1}^{N} \ell(X^i_t) .  
\]
Based on the definition of $y_t$ we have,
\[
\begin{aligned}
	    \bar{\health_t}  &= \frac{1}{N}  \sum_{i=1}^{N} \sum_j [\ind\{ X_t^i = (\oplus, j) \} - \ind\{ X_t^i = (\ominus, j) \} ]\\
	   		  &= \frac{1}{N}  \biggl( 2\sum_{i=1}^{N}  \sum_j \ind\{ X_t^i = (\oplus, j) - N \biggr) \\
		          &= 2\tily_t
\end{aligned}
\]

Let $\health_{\tau}^{C_i}$ denote the average of indicators in class $i$, and $\Health_t^{C_i}$ denote the discounted sum of   $\health_{\tau}^{C_i}$. 
\[
  \Health_{t}^{C_j} =  \sum_{k=0}^{\tau} \beta^k \health_{\tau+1-k}^{C_j} 
  \] 
Note that at time $t$ only $\health_{\tau}^{C_j}$ is updated; other  $\health_{\tau}^{C_i}, \; i\neq j$ remain at their values at $t-1$. Therefore we have
 $ \Health_{m\tau+j + n}^{C_j} =\Health_{m\tau+j}^{C_j}$, for all $j$ and $n \in \{1,2,\cdots,m-1\}$. 
Consequently,
\[
\begin{aligned}
	\bar{\health_t}  &= \frac{1}{N} \sum_{i=1}^{N} \ell(X^i_t) \\
		&=  \frac{1}{N} \left[ \frac{N}{m} \left(\sum_{i=j+1}^{m-1}  \health_{\tau}^{C_i}+ \sum_{i=0}^{j} \health_{\tau+1}^{C_i} \right) \right] \\
		&=  \frac{1}{m}  \left( \sum_{i=j+1}^{m-1}  \health_{\tau}^{C_i} + \sum_{i=0}^{j} \health_{\tau+1}^{C_i} \right)
\end{aligned}
\]

Therefore, the average QoS at time $t$ is
\[
\begin{aligned}
	\bar{\Health_t} &= \frac{1}{N}  \sum_{i=1}^{N} \Health_t^i\\
		&=  \frac{1}{N} \left[ \frac{N}{m} \left( \Health_t^{C_1} + \cdots + \Health_t^{C_m}      \right)  \right] \\
	     	&= \frac{1}{m}   \sum_{k=0}^{\tau}\beta^k \left(  \sum_{i=j+1}^{m-1}  \health_{\tau-k}^{C_i} + \sum_{i=0}^{j} \health_{\tau+1-k}^{C_i} \right) \\
		&=  \sum_{k=0}^{\tau} \beta^k \bar{\health}_{m(\tau-k)+j} 
\end{aligned}
\]

If tracking is perfect as assumed in the proposition, then $\tily_t = r_t$, or $\bar{\health_t} = 2 r_t$. We obtain from this and \eqref{e:R^beta} the desired conclusion that
$
\bar{\Health_t} = 2R_t^{\beta}$.
\qed

\end{document}